\documentclass[12pt]{amsart}

\usepackage{amsmath, amscd, graphicx, latexsym, hyperref, times, rlepsf}

\oddsidemargin 0in \theoremstyle{plain} \topmargin 0in

\newtheorem{Thm}{Theorem}[section]
\newtheorem{Lem}[Thm]{Lemma}
\newtheorem{Cor}[Thm]{Corollary}
\newtheorem{Prop}[Thm]{Proposition}
\newtheorem{Def}[Thm]{Definition}
\newtheorem{Rem}[Thm]{Remark}
\newtheorem{Ex}[Thm]{Example}

\newcommand{\bfc}{{\mathbb{C}}}

\newcommand{\OB}{{\mathrm{ob}}}

\numberwithin{equation}{section}

\begin{document}

\title[]{Surgery diagrams for horizontal contact structures}

\author{Burak Ozbagci}

\address{Department of Mathematics \\ Ko\c{c} University \\
Istanbul \\ Turkey}

\email{bozbagci@ku.edu.tr}


\date{\today}


\begin{abstract}

We describe Legendrian surgery diagrams for some horizontal
contact structures on non-positive plumbing trees of oriented
circle bundles over spheres with negative Euler numbers. As an
application we determine Milnor fillable contact structures on
some Milnor fillable $3$-manifolds.

\end{abstract}

\maketitle 

\section{Introduction}
An open book on a plumbing of oriented circle bundles over spheres
according to a tree is called \emph{horizontal} if its binding is
a collection of some fibers in the circle bundles, its pages
(excluding the binding) are positively transverse to the fibers of
the circle bundles and the orientation induced on the binding by
the pages coincides with the orientation of the fibers induced by
the fibration. Similarly, a contact structure on such a plumbing
is called \emph{horizontal} if the contact planes (which are
oriented by the differential of the contact form) are positively
transverse to the fibers of the bundles involved in the plumbing.
We will call a plumbing tree \emph{non-positive} if for every
vertex of the tree the sum of the degree of the vertex and the
Euler number of the bundle corresponding to that vertex is
non-positive. Note that a non-positive plumbing is called a
``plumbing with no bad vertex" in \cite{os}. We will also assume
that the Euler number of a circle bundle in a non-positive
plumbing tree is less than or equal to $-2$.

In \cite{eo}, Etg\"{u} and the author constructed horizontal
planar open books on non-positive plumbing trees of oriented
circle bundles over spheres. (See also Gay's paper \cite{ga}). Our
construction of the open books were very explicit as we clearly
described the bindings, the pages and vector fields whose flows
induced the monodromy maps of the open books. We also showed that
the contact structures compatible with these open books are
horizontal as well. Note that for the existence of compatible
contact structures we had to rely on a result of Thurston and
Winkelnkemper \cite{tw}. Combining with a supplementary result of
Etnyre in \cite{et2} we obtained a somewhat satisfactory
description of these horizontal contact structures. Our goal in
this paper is to produce yet another (perhaps more explicit)
description of these contact structures, up to isomorphism, by
finding Legendrian surgery diagrams for them.

In \cite{sc}, Sch\"{o}nenberger also gives a construction of (not
necessarily horizontal) planar open books compatible with some
Stein fillable contact structures on non-positive plumbing trees.
He starts with the usual surgery diagram of a given plumbing tree
in $S^3$ and then modifies this diagram by some simple handle
slides to put the surgery link in a certain form so that the
components of the new surgery link can be viewed as embedded
curves in distinct pages of some open book in $S^3$. By Legendrian
realizing these curves on the pages and applying Legendrian
surgery, Sch\"{o}nenberger obtains Stein fillable contact
structures along with the open books compatible with these contact
structures. The advantage of his method is that one can identify
the contact structures by their surgery diagrams. It does not,
however, seem possible to tell which one of these contact
structures (for a fixed plumbing tree) is horizontal applying his
method. Our strategy here is to compare the open books obtained by
these two different approaches to determine the Legendrian surgery
diagrams for the horizontal contact structures given in \cite{eo}
on non-positive plumbing trees.

A contact $3$-manifold $(Y, \xi)$ is said to be \emph{Milnor
fillable} if it is isomorphic to the contact boundary of an
isolated complex surface singularity $(X, x)$. It is a well-known
result of Grauert \cite{gr} that an oriented $3$-manifold has a
Milnor fillable contact structure if and only if it can be
obtained by plumbing oriented circle bundles over surfaces
according to a graph with negative definite intersection matrix.
On the other hand, a recent discovery of Caubel, N\'{e}methi and
Popescu-Pampu (cf. \cite{cnp}) is that any $3$-manifold admits at
most one Milnor fillable contact structure, up to isomorphism. As
an application of our constructions we will be able to determine
this unique Milnor fillable contact structure for some Milnor
fillable $3$-manifolds.

\vspace{1ex}

\noindent {\bf {Acknowledgement}}: We would like to thank
Andr\'{a}s Stipsicz for his valuable comments on the first draft
of this paper. We would also like to thank Tolga Etg\"{u} for
helpful conversations.

\section{Open book decompositions}\label{openbook}

Suppose that for an oriented link $B$ in a closed and oriented
$3$-manifold $Y$ the complement $Y\setminus B$ fibers over the
circle as $\pi \colon Y \setminus B \to S^1$ such that
$\pi^{-1}(\theta) = \Sigma_\theta $ is the interior of a compact
surface with $\partial \Sigma_\theta =B$, for all $\theta \in
S^1$. Then $(B, \pi)$ is called an \emph{open book decomposition}
(or just an \emph{open book}) of $Y$. For each $\theta \in S^1$,
the surface $\Sigma_\theta$ is called a \emph{page}, while $B$ the
\emph{binding} of the open book. The monodromy of the fibration
$\pi$ is defined as the diffeomorphism of a fixed page which is
given by the first return map of a flow that is transverse to the
pages and meridional near the binding. The isotopy class of this
diffeomorphism is independent of the chosen flow and we will refer
to that as the \emph{monodromy} of the open book decomposition. An
open book $(B, \pi)$ on a $3$-manifold $Y$ is said to be
\emph{isomorphic} to an open book $(B^\prime, \pi^\prime)$ on a
$3$-manifold $Y^\prime$, if there is a diffeomorphism $f: (Y,B)
\to (Y^\prime, B^\prime)$ such that $\pi^\prime \circ f = \pi$ on
$Y \setminus B$.

An open book can also be described as follows. First consider the
mapping torus $$\Sigma_\phi= [0,1]\times \Sigma/(1,x)\sim (0,
\phi(x))$$ where $\Sigma$ is a compact oriented surface with $r$
boundary components and $\phi$ is an element of
 the mapping class group $\Gamma_\Sigma$ of $\Sigma$.
 Since
$\phi$ is the identity map on $\partial \Sigma$, the boundary
$\partial \Sigma_\phi$ of the mapping torus $\Sigma_\phi$ can be
canonically identified with $r$ copies of $T^2$. Now we glue in $r$
copies of $D^2\times S^1$ to cap off $\Sigma_\phi$ so that $\partial
D^2$ is identified with $S^1 = [0,1] / (0\sim 1)$ and the $S^1$
factor in $D^2 \times S^1$ is identified with a boundary component
of $\partial \Sigma$. Thus we get a closed $3$-manifold $Y=
\Sigma_\phi \cup_{r} D^2 \times S^1 $ equipped with an open book
decomposition whose binding is the union of the core circles $D^2
\times S^1$'s that we glue to $\Sigma_\phi$ to obtain $Y$. In
conclusion, an element $\phi \in \Gamma_\Sigma$ determines a
$3$-manifold together with an ``abstract" open book decomposition on
it. Notice that by conjugating the monodromy $\phi$ of an open book
on a $3$-manifold $Y$ by an element in $\Gamma_\Sigma$ we get an
isomorphic open book on a $3$-manifold $Y^\prime$ which is
diffeomorphic to $Y$.

Suppose that an open book decomposition with page $\Sigma$ is
specified by $\phi \in \Gamma_\Sigma$. Attach a $1$-handle to the
surface $\Sigma$ connecting two points on $\partial \Sigma$ to
obtain a new surface $\Sigma^{\prime}$. Let $\gamma$ be a closed
curve in $\Sigma^{\prime}$ going over the new $1$-handle exactly
once. Define a new open book decomposition with $ \phi^\prime=
\phi \circ t_\gamma \in \Gamma_{\Sigma^{\prime}} $, where
$t_\gamma$ denotes the right-handed Dehn twist along $\gamma$. The
resulting open book decomposition is called a \emph{positive
stabilization} of the one defined by $\phi$. Notice that although
the resulting monodromy depends on the chosen curve $\gamma$, the
$3$-manifold specified by $(\Sigma^\prime, \phi^\prime)$ is
diffeomorphic to the $3$-manifold specified by $(\Sigma, \phi)$.

Our interest in finding open books on $3$-manifolds arises from
their connection to contact structures, which we will recall here
very briefly. We will assume throughout this paper that a contact
structure $\xi=\ker \alpha$ is coorientable (i.e., $\alpha$ is a
global $1$-form) and positive (i.e., $\alpha \wedge d\alpha >0 $ ).

{\Def \label{compatible} An open book decomposition $(B,\pi)$ of a
$3$-manifold $Y$ and a contact structure $\xi$ on $Y$ are called
\emph{compatible} if $\xi$ can be represented by a contact form
$\alpha$ such that $\alpha ( B) > 0$ and $d \alpha > 0$  on every
page.}

\vspace{1ex}

\noindent Thurston and Winkelnkemper \cite{tw} showed that every
open book admits a compatible contact structure and conversely
Giroux \cite{gi} proved that every contact $3$-manifold admits a
compatible open book. We refer the reader to \cite{et2} and
\cite{ozst} for more on the correspondence between open books and
contact structures.

\begin{figure}[ht]
  \relabelbox \small {
  \centerline{\epsfbox{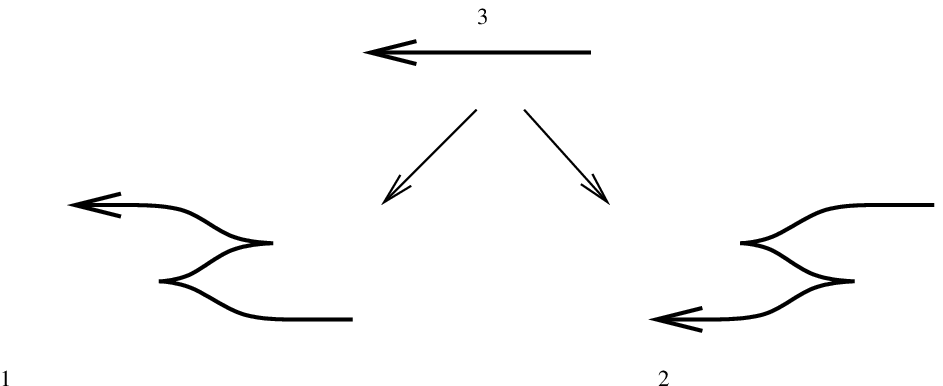}}}

  \relabel{1}{negative stabilization of $K$}

  \relabel{2}{positive stabilization of $K$}
  \relabel{3}{$K$}
  \endrelabelbox
        \caption{} \label{zig}
 \end{figure}

Let $K$ be a Legendrian knot in the standard contact $S^3$ given
by its front projection. We define the positive and negative
stabilization of $K$ as follows: First we orient the knot $K$ and
then if we replace a strand of the knot by an up (down, resp.)
cusp by adding a zigzag as in Figure~\ref{zig} we call the
resulting Legendrian knot the negative (positive, resp.)
stabilization of $K$. Notice that stabilization is a well defined
operation, i.e., it does depend at what point the stabilization is
done.

\begin{figure}[ht]
  \relabelbox \small {
  \centerline{\epsfbox{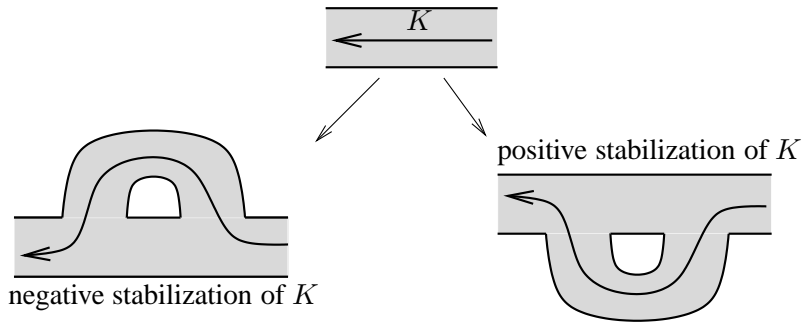}}}

  \relabel{1}{negative stabilization of $K$}

  \relabel{2}{positive stabilization of $K$}
  \relabel{3}{$K$}
  \endrelabelbox
        \caption{Stabilization of a page to include the stabilization of the
   Legendrian knot $K$} \label{legenstab}
 \end{figure}

{\Lem (\cite{et1}) \label{stabb} Let $K$ be a Legendrian knot in a
page of an open book $\OB$ in $S^3$ compatible with its standard
tight contact structure. Then the stabilization of $K$ lies in a
page of another open book in $S^3$ obtained by stabilizing $\OB$.}

In fact, the stabilization of $K$ is obtained by sliding $L$ over
the 1-handle (which is attached to the page to stabilize $\OB$).
Notice that there is a positive and a negative stabilization of
the oriented Legendrian knot $K$ defined by adding a down or an up
cusp, and this choice corresponds to adding a \emph{left} (i.e.,
to the left-hand side of the oriented curve $K$) or a \emph{right}
(i.e., to the right-hand side of the oriented curve $K$) 1-handle
to the surface respectively as shown in Figure~\ref{legenstab}.

\section{Surgery diagrams for horizontal contact
structures} \label{diagrams}

\subsection{Circle bundles with negative Euler numbers}\label{neg}

Let $Y_n$ denote the oriented circle bundle over $S^2$ with Euler
number $e(Y_n)= n < -1$. Then the page of the horizontal open book
$\OB_n$ on $Y_n$ explicitly constructed in \cite{eo} is a planar
surface with $\vert n \vert$ boundary components and the monodromy
is a product of $\vert n \vert $ right-handed Dehn twists, each of
which is along a curve parallel to a boundary component. Now
consider the contact structure $\xi_n$ compatible with the planar
horizontal open book $\OB_n$. In \cite{eo} it was shown that
$\xi_n$ is horizontal as well.

 {\Lem The contact structures
$\xi_n$ is given, up to isomorphism, by the Legendrian surgery
diagram depicted in Figure~\ref{planaropenb}.}

\begin{proof} We will show that the open book compatible with the
contact structure depicted in Figure~\ref{planaropenb} is
isomorphic to $\OB_n$. Consider the core circle $C$ of the open
book $\OB_H$ in $S^3$ given by the positive Hopf link $H$. The
page of $\OB_H$ is an annulus and its monodromy is a right-handed
Dehn twist along $C$. First we Legendrian realize $\gamma$ on a
page of $\OB_H$. Then we use Lemma~\ref{stabb}, to stabilize
$\gamma$ so that the stabilized knot will still be embedded on a
page of another open book in $S^3$. Note that there are two
distinct ways of stabilizing an open book corresponding to two
distinct ways of stabilizing a Legendrian knot in the standard
contact $S^3$. We apply this trick to stabilize $\OB_H$, $ (\vert
n\vert-1 )$-times, by successively attaching $1$-handles while
keeping the genus of the page to be zero.  Here we choose to
stabilize $\OB_H$ always on the same side of the core circle.
\begin{figure}[ht]
  \relabelbox \small {
  \centerline{\epsfbox{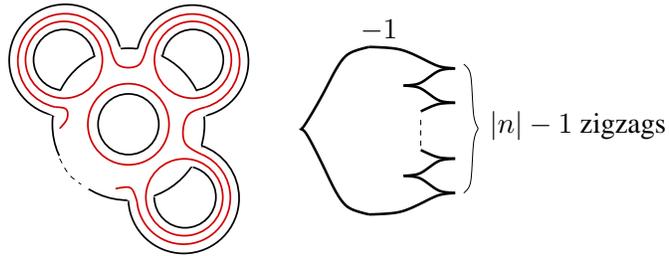}}}

  \relabel{2}{{$\vert n \vert -1$ zigzags}}

  \relabel{4}{{$-1$}}

  \endrelabelbox
        \caption{A page of a planar open book on the left compatible with
        the contact structure on the right.
        The monodromy is given by the product of right-handed Dehn
        twists along the thicker (red) curves.}

        \label{planaropenb}
 \end{figure}

Then we slide $C$ over all the attached $1$-handles, Legendrian
realize the resulting knot on the page of the stabilized open book
and perform Legendrian surgery on this knot to get an open book
compatible with the contact structure induced by the Legendrian
surgery diagram in Figure~\ref{planaropenb}. We observe that the
resulting compatible open book is isomorphic to $\OB_n$ since its
page is a sphere with $\vert n \vert$ holes and its monodromy is a
product of one right-handed Dehn twist along each boundary
component.

\end{proof}

 Note that by reversing
the orientations of the pages (and hence the orientation of the
binding) of $\OB_n$ we get another planar open book
$\overline{\OB}_n$ on $Y_n$. The open book $\overline{\OB}_n$ is
in fact isomorphic to $\OB_n$, since $\OB_n$ and
$\overline{\OB}_n$ have identical pages and the same monodromy map
measured with the respective orientations. To see that
$\overline{\OB}_n$ is also horizontal we simply reverse the
orientation of the fiber (to agree with the orientation of the
binding) as well as the orientation of base $S^2$ of the circle
bundle $Y_n$, so that we do not change the orientation of $Y_n$.
In addition we observe that the contact structure
$\overline{\xi}_n$ compatible with $\overline{\OB_n}$ can be
obtained from $\xi_n$ by changing the orientations of the contact
planes since $\overline{\OB}_n$ is obtained from $\OB_n$ by
changing the orientations of the pages. This immediately implies
that $c_1(\overline{\xi}_n)=-c_1(\xi_n)$. We conclude that the
contact structure induced by the mirror image of the Legendrian
surgery diagram in Figure~\ref{planaropenb} is isomorphic to
$\overline{\xi_n}$. Finally note that $\xi_n$ is isomorphic but
not homotopic to $\overline{\xi}_n$.


\subsection{Non-positive plumbing trees of circle bundles}

In \cite{eo}, an explicit planar horizontal open book $\OB_\Gamma$
on a given non-positive plumbing tree $\Gamma$ of oriented circle
bundles over spheres was constructed. Moreover it was shown that the
compatible contact structure $\xi_\Gamma$ is horizontal as well. In
fact as in the previous section by changing the orientations of the
pages of this horizontal open book and the orientations of the
corresponding contact planes we get another horizontal contact
structure $ \overline{\xi}_\Gamma$. In the following, we will find
Legendrian surgery diagrams for these horizontal contact structures.
We will first consider the easier case of lens spaces and then give
a construction for the general case.

\subsubsection{Lens spaces}\label{le}
The lens space $L(p,q)$ for $p>q \geq 1$ is defined by a
$(-\frac{p}{q})$--surgery on an unknot in $S^3$. Equivalently
$L(p,q)$ can be obtained by a linear plumbing $\Gamma$ of oriented
circle bundles over spheres with Euler numbers $n_1, n_2,\ldots,
n_k$ (see Figure~\ref{rolledup}) where $[n_1, n_2,\ldots, n_k]$ is
the continued fraction expansion for $(-\frac{p}{q})$ with $n_i
\leq -2$, for all $i=1, \cdots, k$.

\begin{figure}[ht]
  \relabelbox \small {
  \centerline{\epsfbox{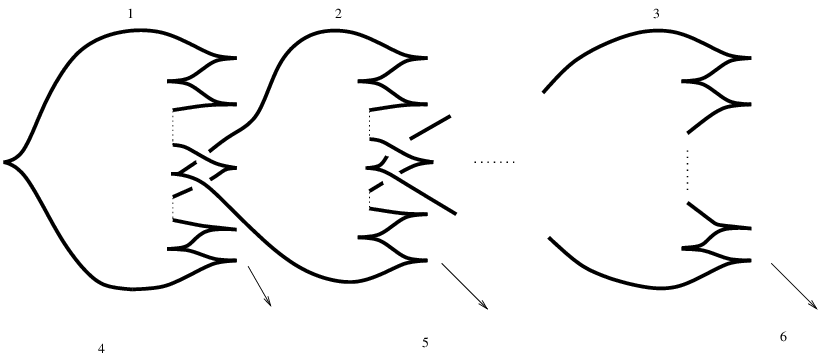}}}
  \relabel{1}{{$-1$}}
  \relabel{2}{{$-1$}}
  \relabel{3}{{$-1$}}
  \relabel{4}{{$\vert n_1 \vert -1$ zigzags}}
  \relabel{5}{{$\vert n_2 \vert -1$ zigzags}}
  \relabel{6}{{$\vert n_k \vert -1$ zigzags}}

  \endrelabelbox
        \caption{}
        \label{lens-horizontal}
\end{figure}

{\Prop The horizontal contact structure $\xi_\Gamma$ is given, up
to isomorphism, by the Legendrian surgery diagram depicted in
Figure~\ref{lens-horizontal}.}

\begin{proof}
We apply the recipe in \cite{eo} to construct a horizontal open
book $\OB_\Gamma$ on $L(p,q)$, with respect to its plumbing
description. First we consider the horizontal open book $\OB_n$
(see Figure~\ref{planaropenb}) we constructed on a circle bundle
with a negative Euler number $n$. Then we glue the open books
$\OB_{n_i}$ corresponding to the circle bundles represented in the
linear plumbing presentation of $L(p,q)$. The key point is that
when we glue two boundary components with boundary parallel
right-handed Dehn twists we end up with only one right-handed Dehn
twist on the connecting neck, as illustrated abstractly in
Figure~\ref{horizopen}. The monodromy of the open book
$\OB_\Gamma$ in Figure~\ref{horizopen} is given by the product of
right-handed Dehn twists along the thicker (red) curves.

\begin{figure}[ht]
  \relabelbox \small {
  \centerline{\epsfbox{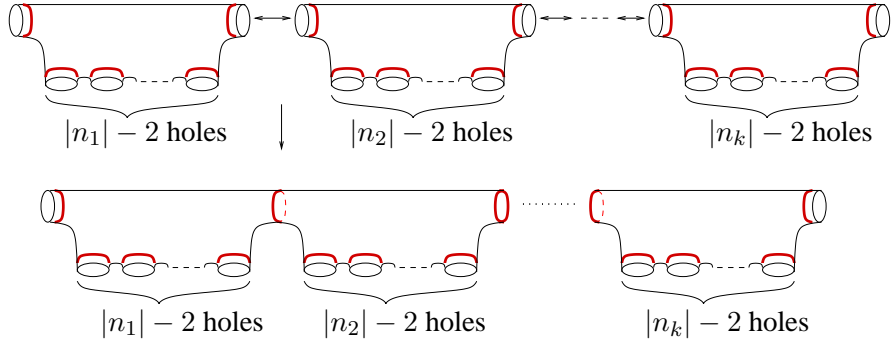}}}
  \relabel{1}{{$\vert n_1 \vert -2$ holes}}
  \relabel{2}{{$\vert n_2 \vert -2$ holes}}
  \relabel{3}{{$\vert n_k \vert -2$ holes}}
  \relabel{4}{{$\vert n_1 \vert -2$ holes}}
  \relabel{5}{{$\vert n_2 \vert -2$ holes}}
  \relabel{6}{{$\vert n_k \vert -2$ holes}}

  \endrelabelbox
        \caption{We illustrate how to plumb open books on circle bundles to get
        the
        horizontal open book $\OB_\Gamma$ in $L(p,q)$. }
        \label{horizopen}
\end{figure}

Now from \cite {sc} we recall how to ``roll up'' a linear plumbing
tree $\Gamma.$ Let $\Gamma$ be the linear plumbing tree with
vertices $Y_{n_1},\ldots, Y_{n_k}$ where each $Y_{n_i}$ is plumbed
only to $Y_{n_{i-1}}$ and $Y_{n_{i+1}}$, for $i=2,\ldots, k-1$, as
shown on the left-hand side of Figure~\ref{rolledup}. The standard
surgery diagram for $\Gamma$ is a chain of unknots $U_1, \ldots,
U_k$ with each $U_i$ simply linking $U_{i-1}$ and $U_{i+1},
i=2,\ldots, k-1$ such that $U_i$ has framing $n_i$. We think of
this chain as horizontal with components labelled from left to
right. Let $U_1'=U_1.$ Slide $U_2$ over $U_1$ to get a new link
with $U_2$ replaced by an unknot $U_2'$ that now links $U_1,$
$n_1+1$ times. Now slide $U_3$ over $U_2'.$ Continue in this way
until $U_k$ is slid over $U_{k-1}'.$ The new link $L$ is called
the ``rolled up" surgery diagram as depicted on the right-hand
side of Figure~\ref{rolledup}.
\begin{figure}[ht]
  \relabelbox \small {
  \centerline{\epsfbox{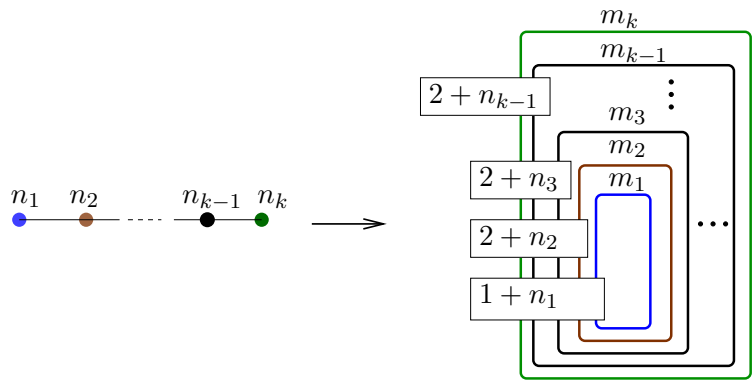}}}
  \relabel{1}{{$n_1 $}}
  \relabel{2}{{$n_2 $}}
  \relabel{3}{{$ n_{k-1} $}}
  \relabel{4}{{$n_k $ }}
  \relabel{5}{{$m_1$ }}
  \relabel{6}{{$m_2$ }}
 \relabel{7}{{$m_3$ }}
 \relabel{8}{{$m_{k-1}$ }}
 \relabel{9}{{$m_k$ }}
 \relabel{a}{{$1+n_1$ }}
 \relabel{b}{{$2+n_2$ }}
 \relabel{c}{{$2+n_3$ }}
 \relabel{d}{{$2+n_{k-1}$ }}
  \endrelabelbox
        \caption{Linear diagram for $L(p,q)$ on the left, rolled up version on the right.
        The number inside a box indicates the number of full-twists that should be applied
        to the strands entering into that box.}
        \label{rolledup}
\end{figure}

Note that a $2$-handle slide corresponds to a change of basis in
the second homology of the $4$-manifold bounded by the
$3$-manifold at hand. We observe a few simple features of this
construction:

\vspace{1ex}

$\bullet$ For fixed $i$,  $U_j'$ links $U_i'$ the same number of
times for any $j>i$ . If we denote this linking number by $l_i$,
then we have $l_i=n_1+\ldots+n_i +2i-1$.

\vspace{1ex}

$\bullet$ Since $m_{i+1}-m_i=n_{i+1}+2,$ we have $m_i= 2(i-1) +
\Sigma_{j=1}^i n_j$. Note that the framings $m_i$ on the $U_i'$'s
are non-increasing and decrease only when $n_i<-2.$

\vspace{1ex}

$\bullet$ The meridian $\mu_i$ for $U_i$ simply links $U'_i\cup
\ldots \cup U'_k.$

\vspace{1ex}

$\bullet$ $L$ sits in an unknotted solid torus neighborhood of
$U_1.$

Note that we can realize $L$ as a Legendrian link in the standard
contact $S^3$ such that $U'_i$ is the Legendrian push off of
$U'_{i-1}$ with $|n_{i}+2|$ stabilizations. Using this observation
we can find an open book in $S^3$ such that $U'_1 \cup \ldots \cup
U'_k$ are embedded in distinct pages of this open book as follows:
First find an open book in $S^3$ as in Section~\ref{neg} for which
the innermost knot $U'_1$ in $L$ is embedded on a page. Then
stabilize this open book appropriately so that $U'_2$ is also
embedded on a page. We will get the desired open book by
continuing this way. Note that there are two different ways of
stabilizing our open book at every step of this construction. For
our purposes in this section we will always choose to stabilize in
the same direction. Consequently when we Legendrian realize these
knots on distinct pages of our open book in $S^3$, they will all
be stabilized in the same direction in $S^3$. By performing
Legendrian surgeries on them, we get a planar open book on
$L(p,q)$ compatible with the resulting contact structure shown in
Figure~\ref{planaropen}.

\begin{figure}[ht]
  \relabelbox \small {
  \centerline{\epsfbox{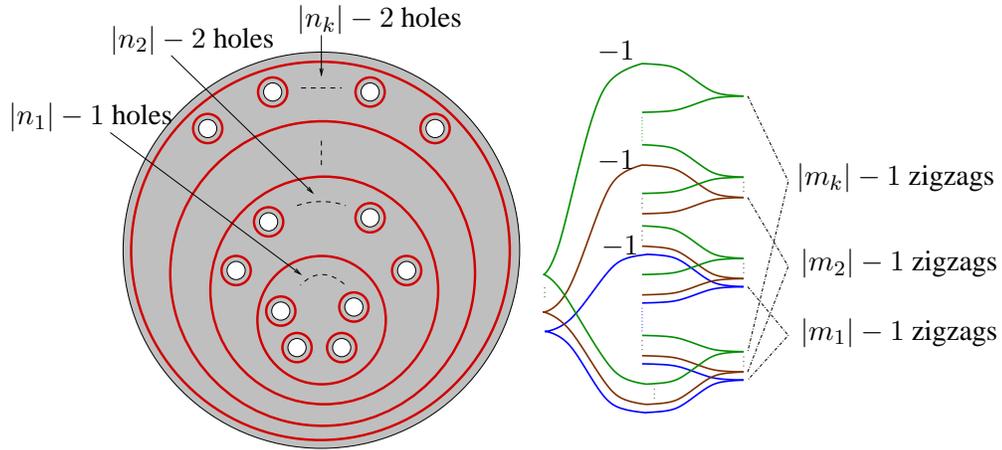}}}
  \relabel{1}{{$\vert n_k \vert -2$ holes}}
  \relabel{2}{{$\vert n_2 \vert -2$ holes}}
  \relabel{3}{{$\vert n_1 \vert -1$ holes}}
   \relabel{4}{{$\vert m_1 \vert -1$ zigzags}}
  \relabel{5}{{$\vert m_2 \vert -1$ zigzags}}
  \relabel{6}{{$\vert m_k \vert -1$ zigzags}}
  \relabel{7}{{$-1$}}
  \relabel{8}{{$-1$}}
  \relabel{9}{{$-1$}}
  \endrelabelbox
        \caption{A page of a planar open book on the left compatible with
        the contact structure on the right.
        The monodromy is given by the product of right-handed Dehn
        twists along the thicker (red) curves.
         }
        \label{planaropen}
\end{figure}

It is not too hard to see that the open book in
Figure~\ref{planaropen} is isomorphic to the horizontal open book
$\OB_\Gamma$ shown in Figure~\ref{horizopen}. Hence we conclude
that the contact structure in Figure~\ref{planaropen} is
isomorphic to the horizontal contact structure $\xi_\Gamma$. Next
we claim that the surgery diagrams in Figure~\ref{lens-horizontal}
and in Figure~\ref{planaropen} induce isomorphic contact
structures. Let us denote these contact structures by $\xi$ and
$\xi'$, respectively.  To prove our claim we first note that these
diagrams induce \emph{diffeomorphic} Stein fillings $(X, J)$ and
$(X', J')$ of the contact $3$-manifolds $(\partial X, \xi)$ and
$(\partial X', \xi')$, where $\partial X \cong \partial X'\cong
L(p,q)$. Now consider $U_i$'s and $U'_i$'s as second homology
classes in these Stein fillings. Then we can identify $U'_i$ as
$U_1+ \ldots +U_i$, by our construction. It is well-known (cf.
\cite{go}) that $$ < c_1(X, J), U_i> =rot(U_i)\; \mbox{and}
<c_1(X', J'), U'_i > =rot(U'_i)$$ where rotation numbers
$rot(U_i)$ and $rot(U'_i)$ are computed using the diagrams in
Figure~\ref{lens-horizontal} and in Figure~\ref{planaropen},
respectively, with appropriate orientations. It follows that we
have
$$< c_1 (X', J'), U_i >=< c_1 (X, J), U_i > $$ for $1 \leq i\leq k$.
Thus the classification of tight contact structures on lens spaces
(cf. \cite{h1}) implies that $\xi$ is isomorphic to $\xi'$.

Similarly, by putting all the zigzags on the left-hand side in
Figure~\ref{lens-horizontal} we get a Legendrian surgery diagram
which induces a contact structure isomorphic to
$\overline{\xi}_\Gamma$. Note that $\xi_\Gamma$ is isomorphic but
not homotopic to $\overline{\xi}_\Gamma$.
\end{proof}

\subsubsection{General case} We would like to explain the main idea in
this section before we get into the details. Given a non-positive
plumbing tree of circle bundles over spheres. We first modify (cf.
\cite{sc}) the usual surgery description of this plumbing tree as
in the previous section to end up with a rolled up version of the
original diagram such that each component of the new link is
embedded in a page of some open book in $S^3$. To find this open
book in $S^3$ we need to appropriately stabilize $\OB_H$ several
times paying attention to the fact that there are two possible
ways of stabilizing at each time. Note that this choice will
determine how to stabilize a Legendrian push-off of a Legendrian
knot $K$ in $S^3$ by adding either a right or a left zigzag (cf.
Lemma~\ref{stabb}) when we slide $K$ over the attached $1$-handle.
Moreover this will also determine the monodromy curves of the open
book. The key observation is that the monodromy of the horizontal
open book $\OB_\Gamma$ constructed in \cite{eo} for a non-positive
plumbing tree $\Gamma$ is given by a product of right-handed Dehn
twists along \emph{disjoint} curves. Thus we will apply the same
strategy as in the previous section so that our choices in the
stabilizations are dictated by this ``disjointness" condition.
Below we give the details of this discussion.

\begin{figure}[ht]
  \relabelbox \small {
  \centerline{\epsfbox{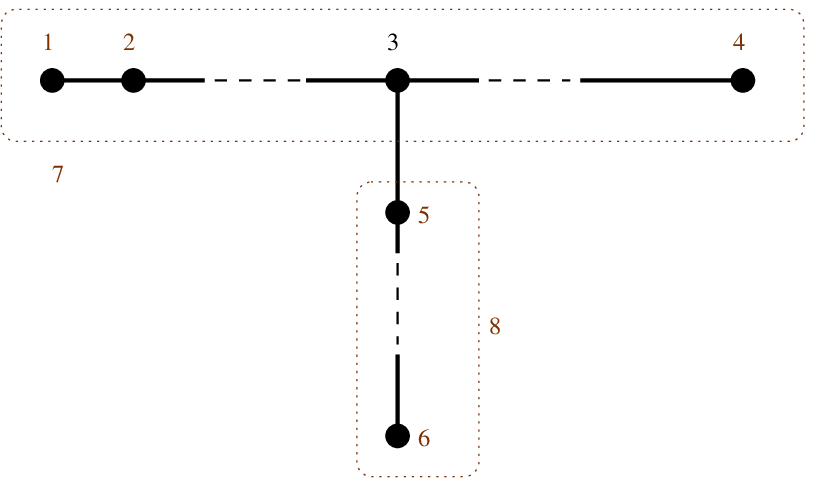}}}
  \relabel{8}{{$\Gamma_2$}}
  \relabel{7}{{$\Gamma_1$}}
  \relabel{3}{{$U_i$}}
 \relabel{1}{{$U_1$}}
 \relabel{2}{{$U_2$}}
 \relabel{4}{{$U_k$}}
 \relabel{5}{{$V_1$}}
 \relabel{6}{{$V_{k'}'$}}
  \endrelabelbox
        \caption{}
        \label{treee}
\end{figure}

Consider a (connected) tree $\Gamma$ which has one degree three
vertex and all the other vertices have degree less than or equal
to two. Thus we can decompose $\Gamma$ into two linear trees
$\Gamma_1$ and $\Gamma_2,$ where the first sphere bundle of
$\Gamma_2$ is plumbed into the $i$th sphere bundle of $\Gamma_1$
as shown in Figure~\ref{treee}.

\begin{figure}[ht]
  \relabelbox \small {
  \centerline{\epsfbox{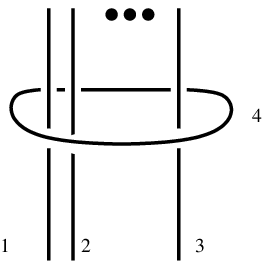}}}
  \relabel{1}{{$U_i'$}}
  \relabel{2}{{$U_{i+1}'$}}
 \relabel{3}{{$U_k'$}}
\relabel{4}{{$V_1'$}}
  \endrelabelbox
        \caption{}
        \label{wrap}
\end{figure}

Let $L_1=U_1'\cup \ldots \cup U_k'$ be the rolled up surgery
diagram of $\Gamma_1= U_1\cup \ldots \cup U_k$. Note that a
neighborhood of the meridian $\mu_i$ for $U_i$ is wrapped once
around $U_i', U_{i+1}', \ldots, U_k'$. Now we identify $V_1$ with
$\mu_i$ and roll up $\Gamma_2$ to obtain $L_2=V_1'\cup \ldots \cup
V'_{k'}$, where $V_1'=V_1$. Here since $V_1$ is identified with
$\mu_i$, $V_1'$ is wrapped once around $U_i', U_{i+1}', \ldots,
U_k'$ as shown in Figure~\ref{wrap}. We will call the resulting
surgery link $L=L_1 \cup L_2$ as the rolled up diagram of
$\Gamma.$

\begin{figure}[ht]
  \relabelbox \small {
  \centerline{\epsfbox{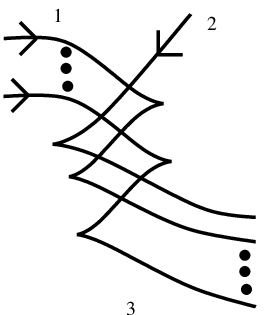}}}
  \relabel{3}{{$U_i'$}}
  \relabel{2}{{$V_1'$}}
 \relabel{1}{{$U_k'$}}

  \endrelabelbox
        \caption{}
        \label{link}
\end{figure}

Now we claim that we can realize $L=L_1 \cup L_2$ as a Legendrian
link using the Legendrian realizations of $L_1$ and $L_2$ we
discussed in the previous section. By the construction there will
be a zigzag in the stabilization of $U_i'$ (and in its subsequent
Legendrian push-offs $U_{i+1}', \ldots, U_k'$ with additional
zigzags) so that we may link $V_1'$ (and hence all of $L_2$) into
$U_i', U_{i+1}', \ldots, U_k'$ using this zigzag as shown in
Figure~\ref{link}.


\begin{figure}[ht]
  \relabelbox \small {
  \centerline{\epsfbox{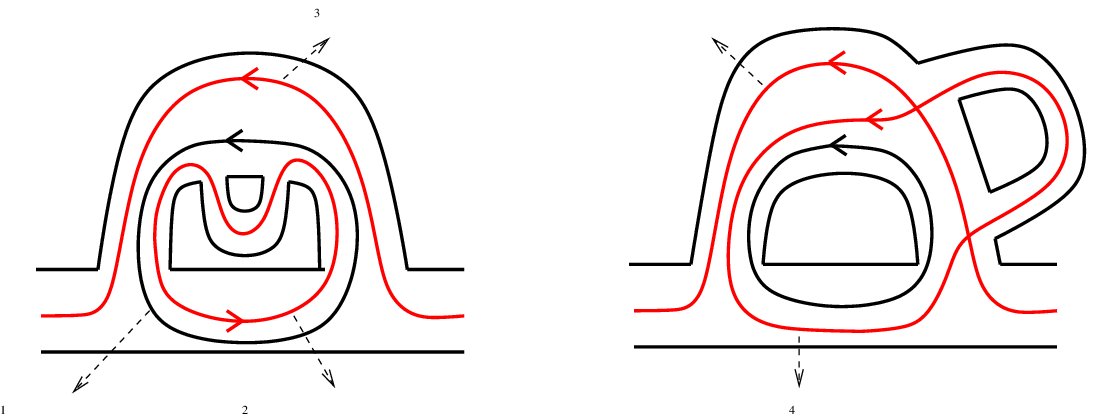}}}
 \relabel{1}{core circle $C$}
 \relabel{2}{once positively stabilized $C$}
 \relabel{3}{once negatively stabilized $K$}
 \relabel{4}{once negatively stabilized $C$}
  \endrelabelbox
        \caption{}
        \label{oppos}
\end{figure}

Hence we can Legendrian realize $L_1$ and $L_2$ as in the previous
section to get a Legendrian realization of $L$. On the other hand,
once we have the rolled up diagram of $\Gamma$ then we can find an
open book in $S^3$ which includes the surgery curves in $L$ on its
pages as follows: Consider the open book corresponding to
$\Gamma_1$ constructed as in Section~\ref{le}. Let us call the
annuli cut out by the large concentric circles in
Figure~\ref{planaropen} as \emph{levels} of this open book.
Observe that the $i$th sphere bundle of $\Gamma_1$ corresponds to
the $i$th level in the open book corresponding to $\Gamma_1$.
Since $V_1$ is linked once to $U_i', U_{i+1}', \ldots, U_k'$ we
can take one of the annulus (a neighborhood of one of the
punctures) in this level as the starting annulus while building up
the open book which includes the surgery curves in $L_2$. We
stabilize  this annulus in the $i$th level \emph{towards the
inside direction} (as many times as necessary) so that we can
embed the curves of $\Gamma_2$ (starting with $V_1$) into the
resulting open book. Note that we need to stabilize towards the
inside direction to keep the curves in $L_1$ disjoint from the
curves in $L_2$. We illustrate this in Figure~\ref{oppos}. On the
left a knot $K$ is stabilized once negatively and then the core
circle $C$ of the 1-handle we use for this stabilization is
stabilized once positively. The two resulting stabilized knots are
clearly disjoint. On the right a knot $K$ is stabilized once
negatively and then the core circle $C$ of the 1-handle we use for
this stabilization is also stabilized negatively once. The two
resulting stabilized knots are clearly \emph{not} disjoint.

Since we have all the surgery curves in $L$ embedded in the pages
of an open book $\OB$ in $S^3$ (compatible with its standard tight
contact structure) now we claim that the open book obtained by
performing Legendrian surgery on $L$ is isomorphic to
$\OB_\Gamma$. Recall that we obtain the open book $\OB_\Gamma$ by
gluing some pieces (see Figure~\ref{horizopen}) which can be
viewed as nothing but plugging in these pieces when converted into
planar diagrams. In order to prove our claim we simply observe
that the result of all the necessary stabilizations of one of the
punctures in the $i$th level of the open book of $\Gamma_1$
towards the inside direction (to embed all the curves of
$\Gamma_2$) is equivalent to ``plugging in" the planar open book
of $\Gamma_2$ into that puncture.

Summarizing this argument what we proved so far is that the
Legendrian surgery diagram for the horizontal contact structure
$\xi_\Gamma$ can be obtained by Legendrian realizing $L_1$ and
$L_2$ such that the direction we stabilize knots in $L_2$ is
opposite to the direction we stabilize knots in $L_1$. By taking
orientations into account as in Figure~\ref{link}, however,  we
realize that we need to put all the zigzags \emph{on the same
side} of all the Legendrian knots that appear in the diagram. That
is because $U_i'$ has down-cups so that $V_1'$ should have
up-cusps. We would like to illustrate the ideas above in the
following example.

\begin{figure}[ht]
  \relabelbox \small {
  \centerline{\epsfbox{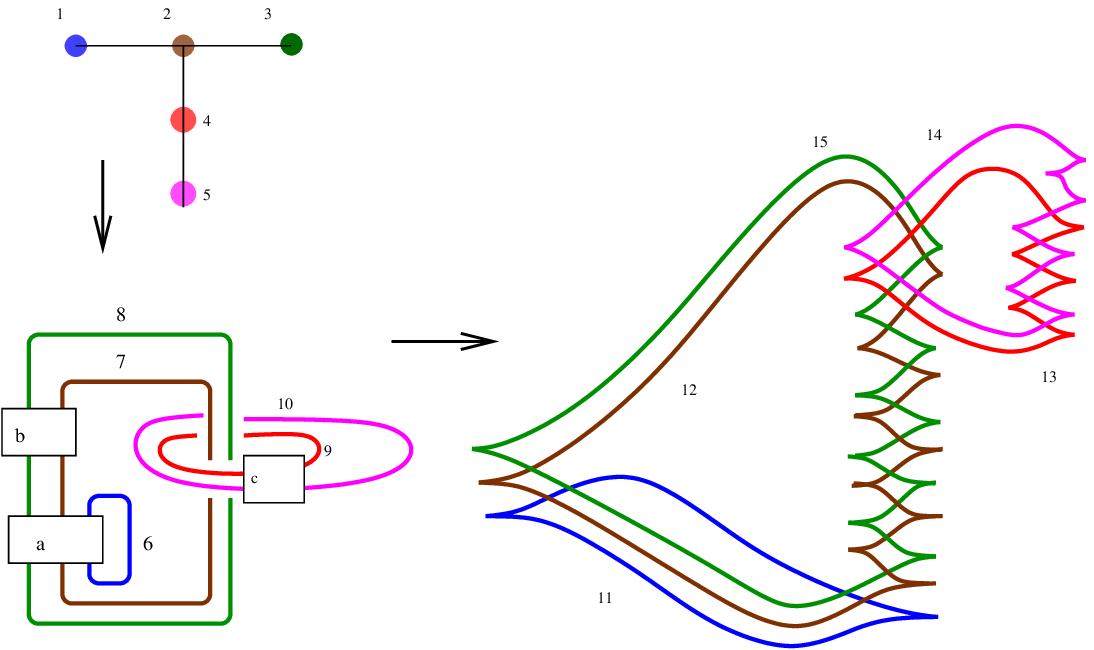}}}
 \relabel{1}{{$-2$}}
  \relabel{2}{{$-6$}}
  \relabel{3}{{$-2$}}
  \relabel{4}{{$-4$}}
  \relabel{5}{{$-3$}}
  \relabel{6}{{$-2$}}
  \relabel{7}{{$-6$}}
  \relabel{8}{{$-6$}}
  \relabel{9}{{$-4$}}
  \relabel{10}{{$-5$}}
  \relabel{a}{{$-1$}}
  \relabel{b}{{$-4$}}
  \relabel{c}{{$-3$}}
 \relabel{11}{{$-1$}}
  \relabel{12}{{$-1$}}
  \relabel{13}{{$-1$}}
  \relabel{14}{{$-1$}}
  \relabel{15}{{$-1$}}
  \endrelabelbox
        \caption{}
        \label{examplerolledup}
\end{figure}

\begin{figure}[ht]
  \relabelbox \small {
  \centerline{\epsfbox{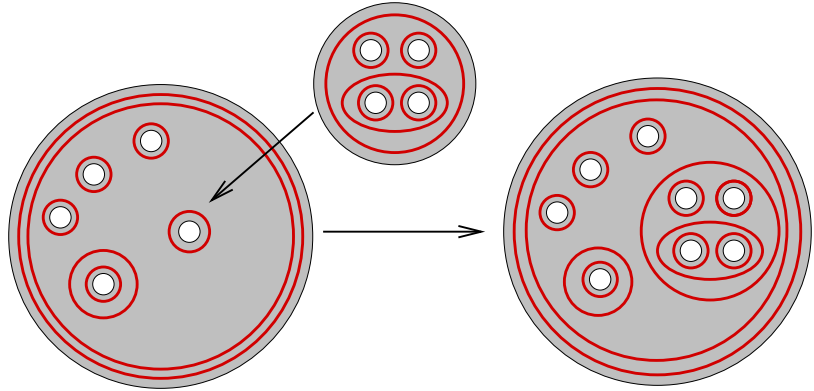}}}

\endrelabelbox
  \caption{On the left: the open book for the linear plumbing of bundles with Euler
numbers $(-2,-6,-2)$; in the middle: the open book for the linear
plumbing of bundles with Euler numbers $(-4,-3)$.
  The horizontal open book on the right-hand side is obtained by the indicated ``plugging in" operation.
 }
        \label{opbook}
\end{figure}

{\Ex \label{roll}  Consider the non-positive plumbing tree and its
rolled up diagram in Figure~\ref{examplerolledup}. On the
right-hand side in Figure~\ref{examplerolledup} we depict the
Legendrian surgery diagram for $\xi_\Gamma$. In
Figure~\ref{opbook} we depict $\OB_\Gamma$. }

\begin{figure}[ht]
  \relabelbox \small {
  \centerline{\epsfbox{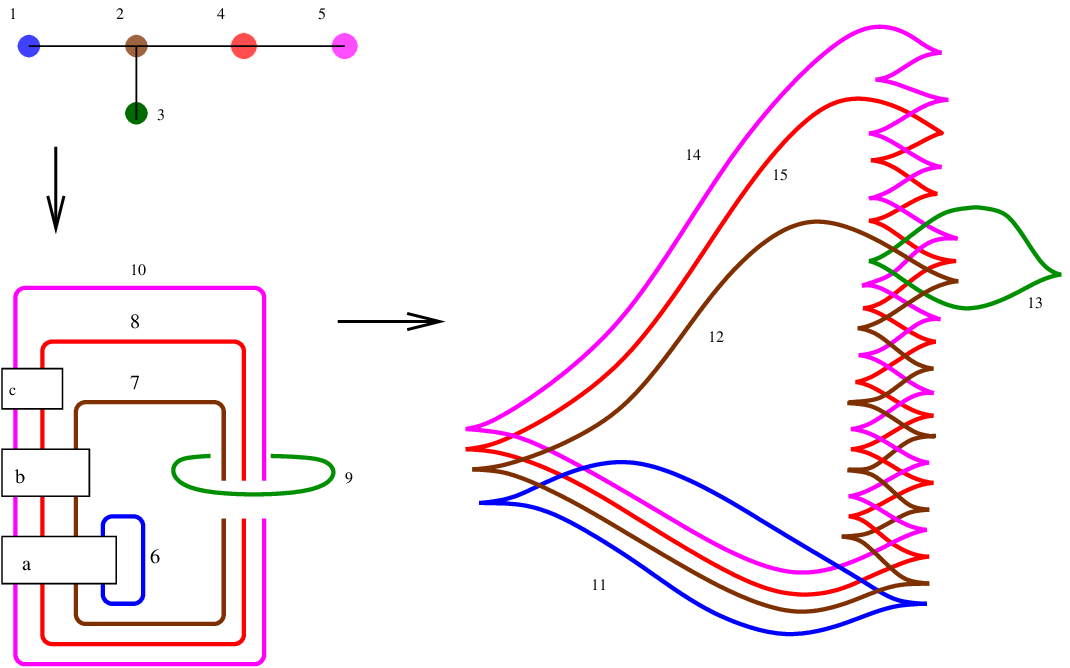}}}
 \relabel{1}{{$-2$}}
  \relabel{2}{{$-6$}}
  \relabel{3}{{$-2$}}
  \relabel{4}{{$-4$}}
  \relabel{5}{{$-3$}}
  \relabel{6}{{$-2$}}
  \relabel{7}{{$-6$}}
  \relabel{8}{{$-8$}}
  \relabel{9}{{$-2$}}
  \relabel{10}{{$-9$}}
  \relabel{a}{{$-1$}}
  \relabel{b}{{$-4$}}
  \relabel{c}{{$-2$}}
 \relabel{11}{{$-1$}}
  \relabel{12}{{$-1$}}
  \relabel{13}{{$-1$}}
  \relabel{14}{{$-1$}}
  \relabel{15}{{$-1$}}
  \endrelabelbox
        \caption{}
        \label{dif}
\end{figure}

{\Rem We may indeed roll up the tree $\Gamma$ in
Example~\ref{roll} in a different way (see Figure~\ref{dif}).
Nevertheless, we end up with a Legendrian surgery diagram which is
isomorphic to the one in Figure~\ref{examplerolledup}, since the
respective corresponding compatible open books are isomorphic.}

So far we analyzed our problem locally, where there is only one
vertex of branching in the plumbing tree. It turns out, however,
the argument can be easily generalized to arbitrary non-positive
plumbing trees. By the non-positivity assumption every time there
is a branching at a vertex there are enough punctures to plug in
the open books corresponding to the linear branches coming out of
that vertex. This argument leads to the following theorem.

{\Thm \label{lege} To obtain a Legendrian surgery diagram for the
horizontal contact structure $\xi_\Gamma$ on a non-positive
plumbing tree $\Gamma$ of circle bundles over spheres we need to
put the rolled up diagram of $\Gamma$ into Legendrian position (as
described above) so that all the zigzags are on the right-hand
side for all the knots in the diagram. }

\section{Milnor fillable contact structures}
Let $(X,x)$ be an isolated complex-analytic singularity. Given a
local embedding of $(X,x)$ in $(\bfc^N, 0)$. Then a small sphere
$S^{2N-1}_{\epsilon} \subset \bfc^N$ centered at the origin
intersects $X$ transversely, and the complex hyperplane distribution
$\xi$ on $Y=X\cap S^{2N-1}_{\epsilon}$ induced by the complex
structure on $X$ is a contact structure. For sufficiently small
radius $\epsilon$ the contact manifold $(Y, \xi)$ is independent of
the embedding and $\epsilon$ up to isomorphism and this isomorphism
type is called the \emph{contact boundary} of $(X,x)$. A contact
manifold $(Y^\prime, \xi^\prime)$ is said to be \emph{Milnor
fillable} and the germ $(X, x)$ is called a \emph{Milnor filling} of
$(Y^\prime, \xi^\prime)$ if $(Y^\prime, \xi^\prime)$ is isomorphic
to the contact boundary $(Y, \xi)$ of $(X, x)$.

It is a well-known result of Grauert \cite{gr} that an oriented
3--manifold has a Milnor fillable contact structure if and only if
it can be obtained by plumbing oriented circle bundles over surfaces
according to a weighted graph with negative definite intersection
matrix. On the other hand, a recent discovery (cf. \cite{cnp}) is
that any 3--manifold admits at most one  Milnor fillable contact
structure, up to isomorphism. Note that any Milnor fillable contact
structure is horizontal (cf. \cite{eo}, Remark 14). The following
proposition was proved in \cite{eo}.

{\Prop  \label{com} Let $Y_\Gamma$ be a $3$-manifold obtained by a
plumbing of circle bundles over spheres according to a tree
$\Gamma$. Suppose that the inequality
$$e_i +2d_i \leq 0$$ holds for every vertex of $\Gamma$, where $e_i$
denotes the Euler number and $d_i$ denotes the degree of the $i$th
vertex. Then $Y_\Gamma$ is Milnor fillable and the planar
horizontal open book $\OB_\Gamma$ is compatible with the unique
Milnor fillable contact structure on $Y_\Gamma$.}

{\Cor A Legendrian surgery diagram for the unique Milnor fillable
contact structure on $Y_\Gamma$ (where $\Gamma$ satisfies the
inequality in Proposition~\ref{com}) can be explicitly determined
using the algorithm presented in Theorem~\ref{lege}.}

\begin{proof}

The horizontal contact structure $\xi_\Gamma$ compatible with
$\OB_\Gamma$ is isomorphic to the unique Milnor fillable contact
structure on $Y_\Gamma$ by Proposition~\ref{com} and $\xi_\Gamma$
can be explicitly determined  using Theorem~\ref{lege} since the
inequality $e_i +2d_i \leq 0$ trivially implies non-positivity for
$\Gamma$.

\end{proof}

{\Ex The Legendrian surgery diagram in
Figure~\ref{examplerolledup} represents the unique Milnor fillable
contact structure, up to isomorphism,  on the given non-positive
plumbing tree. Note that by reversing the way we stabilized the
knots in Figure~\ref{examplerolledup} we end up with an isomorphic
but non-isotopic Milnor fillable contact structure.}







\begin{thebibliography}{99999}

\bibitem{cnp}
C. Caubel, A. N\'{e}methi, P. Popescu-Pampu, \emph{Milnor open books
and Milnor fillable contact 3--manifolds}, Topology {\bf 45} (2006),
no. 3, 673--689.

\bibitem{eo}
T. Etg\"{u} and B. Ozbagci, {\em Explicit horizontal open books on
some plumbings},  Internat. J. Math. {\bf 17} (9) (2006),
1013--1031.

\bibitem{et1}
J. B. Etnyre, {\em Planar open book decompositions and contact
structures,} Internat. Math. Res. Notices {\bf 79} (2004),
4255--4267.

\bibitem{et2}
J. B. Etnyre, {\em Lectures on open book decompositions and
contact structures}, Lecture notes from the Clay Mathematics
Institute Summer School on Floer Homology, Gauge Theory, and Low
Dimensional Topology at the Alfr\'{e}d R\'{e}nyi Institute;
arXiv:math.SG/0409402.

\bibitem{ga}
D. Gay, \emph{Open books and configurations of symplectic
surfaces,} Algebr. Geom. Topol. {\bf3} (2003), 569--586.


\bibitem{gi}
E. Giroux, \emph{G\'{e}ometrie de contact: de la dimension trois
vers les dimensions sup\'{e}rieures,} Proceedings of the
International Congress of Mathematicians (Beijing 2002), Vol. II,
405--414.

\bibitem{go}
R. Gompf, {\em Handlebody construction of Stein surfaces,} Ann. of
Math. {\bf 148} (1998), 619-693.


\bibitem{gr}
H. Grauert, {\em \"Uber Modifikationen und exzeptionelle analytische
Mengen,} Math. Ann. {\bf 146} (1962), 331--368.


\bibitem{h1}
K. Honda, \emph{On the classification of tight contact structures,
I.,} Geom. Topol. {\bf 4} (2000), 309--368 (electronic).



\bibitem{os}
P. Ozsv\'{a}th and Z. Szab\'{o}, \emph{On the Floer homology of
plumbed three-manifolds,} Geom. Topol. {\bf 7} (2003), 185--224
(electronic).

\bibitem{ozst}
B. Ozbagci and A. Stipsicz, {\em Surgery on contact 3--manifolds
and Stein surfaces}, Bolyai Soc. Math. Stud., Vol. {\bf 13},
Springer, 2004.


\bibitem{sc}
S. Sch\"{o}nenberger, {\em Planar open books and symplectic
fillings}, Ph.D. Dissertation, University of Pennsylvania, 2005,
http://www.math.gatech.edu/~etnyre/professionalstuff/gradstudents.html

\bibitem{tw}
W. Thurston and H. Winkelnkemper, \emph{On the existence of contact
forms,} Proc. Amer. Math. Soc. {\bf 52} (1975), 345--347.

\end{thebibliography}
\end{document}